\documentclass{elsart}
\usepackage{amsmath}
\usepackage{amssymb}
\date{}
\journal{}

\newcommand\ce{{\mathcal C}}
\newcommand\Qmax{Q_{\mathrm{max}}}
\newcommand\Qcl{Q_{\mathrm{cl}}}
\newcommand\Qrmax{Q^r_{\mathrm{max}}}

\begin{document}
\begin{frontmatter}

\title{Class of Baer *-rings Defined by a Relaxed Set
of Axioms}

\author{Lia Va\v s}

\address{Department of Mathematics, Physics and Computer Science,
University of the Sciences in Philadelphia, 600 S. 43rd St.,
Philadelphia, PA 19104}

\ead{l.vas@usip.edu}

\begin{abstract}
We consider a class ${\mathcal C}$ of Baer *-rings (also treated
in \cite{Be} and \cite{Lia2}) defined by nine axioms, the last two
of which are particularly strong. We prove that the ninth axiom
follows from the first seven. This gives an affirmative answer to
the question of S. K. Berberian if a Baer *-ring $R$ satisfies the
first seven axioms, is the matrix ring $M_n(R)$ a Baer *-ring.
\end{abstract}


\end{frontmatter}

\section{A Class $\ce$ of Baer *-Rings}
\label{section_class_C}

In \cite{Lia2}, some results on finite von Neumann algebras are
generalized, by purely algebraic proofs, to a certain class $\ce$
of finite Baer *-rings. The the dimension of any module over a
ring from $\ce$ is defined. This dimension is proven to have the
same properties as the dimension for finite von Neumann algebras.
The class $\ce$ is defined via nine axioms, the last two of which
are particulary strong. In this paper, we demonstrate that the
ninth axiom follows from the first seven (Theorem
\ref{A9_is_redundant}). This also gives an affirmative answer to
the question of S. K. Berberian if a Baer *-ring $R$ satisfies the
first seven axioms, is the matrix ring $M_n(R)$ a Baer *-ring.

If $R$ is a Baer *-ring, let us consider the following set of
axioms.

\begin{itemize}
\item[(A1)] $R$ is {\em finite} (i.e. $x^*x=1$ implies $xx^*=1$
for all $x\in R$).

\item[(A2)] $R$ satisfies {\em existence of projections
(EP)-axiom:} for every $0\neq x\in R,$ there exist an self-adjoint
$y\in\{x^*x\}''$ such that $(x^*x)y^2$ is a nonzero projection;

$R$ satisfies the {\em unique positive square root (UPSR)-axiom:}
for every $x\in R$ such that $x=x_1^*x_1+x_2^*x_2+\ldots
+x_n^*x_n$ for some $n$ and some $x_1,x_2,\ldots,x_n\in R$ (such
$x$ is called {\em positive}), there is a unique $y\in\{x^*x\}''$
such that $y^2=x$ and $y$ positive. Such $y$ is denoted by
$x^{1/2}.$

\item[(A3)] Partial isometries are addable.

\item[(A4)] $R$ is {\em symmetric}: for all $x\in R,$ $1+x^*x$ is
invertible.

\item[(A5)] There is a central element $i\in R$ such that $i^2=-1$
and $i^*=-i.$

\item[(A6)] $R$ satisfies the {\em unitary spectral (US)-axiom:}
for each unitary $u\in R$ such that RP$(1-u)=1,$ there exist an
increasingly directed sequence of projections $p_n\in\{u\}''$ with
supremum 1 such that $(1-u)p_n$ is invertible in $p_n Rp_n$ for
every $n.$

\item[(A7)] $R$ satisfies the {\em positive sum (PS)-axiom;} if
$p_n$ is orthogonal sequence of projections with supremum 1 and
$a_n\in R$ such that $0\leq a_n\leq f_n,$ then there is $a\in R$
such that $a p_n=a_n$ for all $n.$
\end{itemize}

S. K. Berberian in \cite{Be} uses these axioms to embed a ring $R$
in a regular ring $Q$. More precisely, he shows the following.

\begin{thm} If $R$ is a Baer *-ring satisfying (A1)-- (A7), then
there is a regular Baer *-ring $Q$ satisfying (A1) -- (A7) such
that $R$ is *-isomorphic to a *-subring of $Q$, all projections,
unitaries and partial isometries of $Q$ are in $R,$ and $Q$ is
unique up to *-isomorphism. \label{RegularRing}
\end{thm}

This result is contained in Theorem 1, p. 217, Theorem 1 and
Corollary 1 p. 220, Corollary 1, p. 221, Theorem 1 and Corollary 1
p. 223, Proposition 3 p. 235, Theorem 1 p. 241, Exercise 4A p. 247
in \cite{Be}.

A ring $Q$ as in Theorem \ref{RegularRing} is called the {\em
regular ring of Baer *-ring $R$.}

In \cite{Lia2}, in addition to Berberian's theorem above, the
following theorem is proven to hold.

\begin{prop} If $R$ is a Baer *-ring satisfying (A1)-- (A7) with $Q$ its regular ring, then
\begin{enumerate}
\item $Q$ is the classical (left and right) ring of quotients
$\Qcl(R)$ of $R.$

\item $Q$ is the maximal (left and right) ring of quotients
$\Qmax(R)$ of $R$ and, thus, self-injective and equal to the (left
and right) injective envelope $E(R)$ of $R.$

\item The ring $M_n(R)$ of $n\times n$ matrices over $R$ is
semihereditary for every positive $n.$
\end{enumerate}
\label{Q=max=classical}
\end{prop}

This result is contained in Proposition 3 and Corollary 5 in
\cite{Lia2}.

If $R$ satisfies (A1) -- (A7), the ring $M_n(R)$ of $n\times n$
matrices is a Rickart *-ring for every $n$ (by Theorem 1, p. 251
in \cite{Be}). In \cite{Be}, Berberian used additional two axioms
in order to ensure that $M_n(R)$ is also a finite Baer *-ring with
the dimension function defined on the set of projections (Theorem
1 and Corollary 2, p. 262 in \cite{Be}). In \cite{Lia2}, Theorem
17 shows that the additional two axioms allow the definition of
dimension to be extended to all the modules over $R$ and that this
dimension has all the nice properties of the dimension studied in
\cite{Lu2} (or \cite{Lu_book}) and \cite{Lia1} for finite von
Neumann algebras.

The additional two axioms are:
\begin{itemize}
\item[(A8)] $M_n(R)$ satisfies the {\em parallelogram law (P):}
for every two projections $p$ and $q,$ \[ p - \inf\{ p,q\}\sim
\sup\{ p,q\}-q.\]

\item[(A9)] Every sequence of orthogonal projections in $M_n(R)$
has a supremum.
\end{itemize}

\begin{defn} Let $\ce$ be the class of Baer *-rings that satisfy
the axioms (A1) -- (A9).
\end{defn}

Every finite $AW^*$-algebra satisfies the axioms (A1) -- (A9)
(remark 1, p. 249 in \cite{Be}). Thus, the class $\ce$ contains
the class of all finite $AW^*$-algebras and, in particular, all
finite von Neumann algebras.

\section{Getting rid of (A9)}
\label{main}

Berberian calls the last two axioms "unwelcome guests" since they
impose conditions on the rings $M_n(R)$ for every $n$ and not just
on the ring $R$ itself. Here we prove that (A9) follows from (A1)
-- (A7), so (A9) is redundant. This gives an affirmative answer to
the question Berberian asked on page 253 (Exercise 4D) in
\cite{Be}: If $R$ satisfies (A1)--(A7), is $M_n(R)$ a Baer *-ring?

\begin{thm} If $R$ satisfies (A1) -- (A7), then $M_n(R)$ is a Baer *-ring for every $n$ (thus, (A9) holds).
\label{A9_is_redundant}
\end{thm}

\begin{pf} Let $R$ be a Baer *-ring satisfying (A1) --
(A7). Then $R$ is semihereditary and the regular ring $Q$ of $R$
is both the classical and the maximal ring of quotients by
Proposition \ref{Q=max=classical}. In this case $Q$ is also the
flat epimorphic hull of $R$ (for detailed review of this notion
see \cite{Stenstrom}) by Example on page 235 in \cite{Stenstrom}.
This fact allows us to use the result (Theorem 2.4, page 54, in
\cite{Evans}) of M. W. Evans:

The following are equivalent for any ring $R$:
\begin{itemize}
\item[i)] $R$ is right semihereditary and the maximal right ring
of quotients $\Qrmax(R)$ is the left and right flat epimorphic
hull of $R$.

\item[ii)] $M_n(R)$ is a right strongly Baer ring for all $n.$
\end{itemize}

Our ring $R$ satisfies the condition i) by discussion above. Since
every strongly Baer ring is Baer (see \cite{Evans}, page 53),
$M_n(R)$ is Baer for every $n$. $M_n(R)$ is a Rickart *-ring
(Theorem 1, p. 251 in \cite{Be}) which is Baer and so it is a Baer
*-ring (Proposition 1, p. 20, \cite{Be}).
\end{pf}

This proves that the axiom (A9) can be avoided and, thus, that it
does not need to be assumed in Chapter 9 of Berberian's book
\cite{Be} and in my paper \cite{Lia2}.

\begin{cor} A Baer *-ring is
in $\ce$ if and only if it satisfies (A1) -- (A8).
\end{cor}

\section{Question}
\label{section_questions}

This still leaves open the question whether (A8), the other
"unwelcome guest", must be assumed. In \cite{Be}, (A8) is used to
show that $M_n(R)$ is finite and in \cite{Lia2}, it is used to
show that the dimension function can be extended from projections
in $R$ to projections in $M_n(R).$ Berberian also asks this
question in \cite{Be} (ch. 9, sec. 57, p. 254). The author shares
the view of Berberian who believes that it might not be necessary
to assume (A8) (p. 254, \cite{Be}), but do not have a proof that
(A8) follows from (A1) -- (A7) at this point.

\end{document}